\documentclass[11pt,fceqn]{article}
\usepackage{amsmath}
\usepackage{color}
\usepackage{cite}
\usepackage{amsfonts}
\usepackage{amssymb,graphics,graphicx,subfigure,caption2,rotating}
\usepackage{amsmath, latexsym}
\usepackage[amsmath,thmmarks]{ntheorem}
\usepackage{setspace} 
\numberwithin{equation}{section}
\newtheorem{theorem}{Theorem}[section]
\newtheorem{lemma}{Lemma}[section]

\newtheorem{cor}{Corollary}[section]

\newtheorem{definition}{Definition}[section]
\newcommand{\beq}{\begin{equation}}
\newcommand{\eeq}{\end{equation}}
\newcommand{\beqn}{\begin{eqnarray}}
\newcommand{\eeqn}{\end{eqnarray}}
\usepackage[margin=1in]{geometry}
\date{}

\begin{document}

\date{}
\title{On the unique solution of the generalized absolute value equation
\thanks{This research was supported by National Natural Science Foundation of
China (No. 11961082).}}
\author{Shi-Liang Wu\thanks{Corresponding author: wushiliang1999@126.com}, Shu-Qian Shen\thanks{sqshen@upc.edu.cn}\\
{\small{\it School of Mathematics, Yunnan Normal University, }}\\
{\small{\it  Kunming, Yunnan, 650500, PR China}}\\
{\small{\it $^{\dagger}$College of Science, China University of Petroleum,  }}\\
{\small{\it Qingdao, Shandong, 266580, PR China}}
}

 \maketitle
\begin{abstract}
In this paper, some useful necessary and sufficient conditions for
the unique solution of the  generalized absolute value equation
(GAVE) $Ax-B|x|=b$ with $A, B\in \mathbb{R}^{n\times n}$ from the
optimization field are first presented, which cover the fundamental
theorem for the unique solution of the linear system $Ax=b$ with
$A\in \mathbb{R}^{n\times n}$. Not only that, some new sufficient
conditions for the unique solution of the GAVE are obtained, which
are weaker than the previous published works.


\end{abstract}

\section{Introduction}
In this paper, we concentrate on  the generalized absolute value
equations (GAVE), whose form is below
\begin{equation}\label{eq:1}
Ax+B|x|=b,
\end{equation}
with $A, B\in \mathbb{R}^{n\times n}$ and $b\in \mathbb{R}^{n}$.
When $B = I$, where $I$ stands for the identity matrix, the GAVE
(\ref{eq:1}) reduces to the absolute value equations (AVE)
\begin{equation}\label{eq:2}
Ax+|x| = b.
\end{equation}
The GAVE (AVE) have received considerable attention because they are
used as a useful tool in the optimization field, such as the linear
complementarity problem, linear programming and convex quadratic
programming, and so on. One can see
\cite{Rohn,Mangasarian,Mangasarian2,Wu,Cottle} for more details.


The research of the unique solution is a very important branch of
theoretical analysis of the GAVE (AVE) because the goal of many
numerical methods is to obtain the unique solution of the GAVE
(AVE), including the modified generalized Newton
method\cite{Li,Lian}, the preconditioned AOR  method \cite{Li2} and
the modified Newton-type method\cite{Wang} and the sign accord
method \cite{Rohn3}. Whereas, by observing the structure of the GAVE
(\ref{eq:1}), it is not difficult to see that the nonlinear and
nondifferentiable term $B|x|$ often leads to the nondeterminacy of
the solution of the GAVE (1.1). In this case, we have to give some
constraints to guarantee that  the GAVE (\ref{eq:1}) has a unique
solution. Recently, some sufficient conditions for the unique
solution of the GAVE (\ref{eq:1}) have been obtained in the
literatures.  For example,  in \cite{Rohn4}, John \emph{et al.
}found that the GAVE (\ref{eq:1}) for any $b\in \mathbb{R}^{n}$ has
a unique solution if $\rho(|A^{-1}B|)<1$, where $\rho(\cdot)$
denotes the spectral radius of the matrix. From  the singular value
of the matrix, John in \cite{Rohn5} showed that the GAVE
(\ref{eq:1}) for any $b\in \mathbb{R}^{n}$ has a unique solution
when $\sigma_{1}(|B|) < \sigma_{n}(A)$, where $\sigma_{1}$ and
$\sigma_{n}$, respectively, denote the maximal and minimal singular
value of the matrix. Based on the work in\cite{Rohn5}, Wu and Li in
\cite{Wu5} obtained an improved result, i.e., if $\sigma_{1}(B) <
\sigma_{n}(A)$, then the GAVE (\ref{eq:1}) for any $b\in
\mathbb{R}^{n}$ has a unique solution.   Other sufficient conditions
for the unique solution of the GAVE (\ref{eq:1}), one can see
\cite{Rohn,Rohn3} for more details.

By investigating the previous published works in
\cite{Rohn5,Rohn4,Rohn,Rohn3,Wu5}, the presented conditions for the
unique solution of the GAVE (\ref{eq:1}) just are sufficient
conditions. Although Wu and Li in \cite{Wu2} presented some
necessary and sufficient conditions for the unique solution of the
AVE (\ref{eq:2}), it is regretful that these results in \cite{Wu2}
are not suitable for the  GAVE (\ref{eq:1}) because matrix $B\in
\mathbb{R}^{n\times n}$ in (\ref{eq:1}) is free. So far, to our
knowledge, the necessary and sufficient condition for the unique
solution of the GAVE (\ref{eq:1}) is \emph{void}, which is our
motivation for this paper.  That is to say, the intent of the
present paper is to address this question,  and present some
necessary and sufficient conditions for the unique solution of the
the GAVE (\ref{eq:1}). Although these results in \cite{Wu2} for the
GAVE (\ref{eq:1}) is invalid, the idea employed in \cite{Wu2} makes
us suffice to establish the necessary and sufficient condition for
the unique solution of the GAVE (\ref{eq:1}). Not only that, these
necessary and sufficient conditions for the unique solution of the
GAVE (\ref{eq:1}) not only contain the fundamental theorem for the
unique solution of the linear system $Ax=b$ with $A\in
\mathbb{R}^{n\times n}$, but also yield some new sufficient
conditions for the unique solution of the GAVE (\ref{eq:1}). These
new sufficient conditions are weaker than the previous published
works.

The rest of the paper unfolds below. Section 2 consists of some
useful lemmas and the definition of $P$-matrix. Section 3 contains
some necessary and sufficient conditions for the unique solution of
the GAVE (\ref{eq:1}). It also contains some new sufficient
conditions for the unique solution of the GAVE (\ref{eq:1}). In
Section 4, some conclusions are given to end the paper.

\section{Preliminaries}
In this section, we introduce the definition of $P$-matrix and some
useful lemmas for the later discussion.

\begin{definition} \emph{\cite{Cottle}}  Matrix $A\in \mathbb{R}^{n\times n}$ is called a $P$-matrix if all
its principal minors are positive.
\end{definition}

\begin{lemma} \emph{\cite{Cottle}}  The linear complementarity problem, which finds  $z\in \mathbb{R}^{n}$ such that
\[
w=Mz+q\geq0,\ z\geq0 \ \mbox{and}\ z^{T}w=0 \ \mbox{with}\ M\in
\mathbb{R}^{n\times n},
\]
 has a unique solution for any $q\in \mathbb{R}^{n}$ if and only if the
matrix $M$ is a $P$-matrix.
\end{lemma}

\begin{lemma} \emph{\cite{Johnson}}
Matrix $M$ is a $P$-matrix if and only if matrix $MD+I-D$ is
nonsingular for any diagonal matrix $D=\mbox{diag}(d_{i})$ with $0
\leq d_{i}\leq1$.
\end{lemma}

By the way, Lemma 2.2 implies that matrix $M$ is a $P$-matrix, which
is equivalent:
\begin{description}
\item (1) matrix $M(I-D)+D$ is
nonsingular for any diagonal matrix $D=\mbox{diag}(d_{i})$ with $0
\leq d_{i}\leq1$.
\item (2)  matrix $MF_{0}+F_{1}$ is
nonsingular, where $F_{0}, F_{1}\in \mathbb{R}^{n\times n}$ are two
arbitrary nonnegative diagonal matrices with $F=F_{0}+F_{1}>0$.
\end{description}

\begin{lemma} \emph{\cite{Loyka}}
Let $A,B\in \mathbb{R}^{n\times n}$. Then
\[\sigma_{i}(A+B)\geq\sigma_{i}(A)-\sigma_{1}(B), i=1,2\ldots,n,
\]
where $\sigma_{1}\geq\ldots\geq\sigma_{n}(\geq 0)$ are the singular
values of matrix.
\end{lemma}

%
%

\section{Main results}
In this section, based on the above results in Section 2, we shall
address the problem of  the necessary and sufficient condition for
the unique solution of the GAVE (\ref{eq:1}).

First, based on Lemma 2.1,  we can obtain Theorem 3.1.

\begin{theorem}
Let $A+B$ be nonsingular in $(\ref{eq:1})$. Then the GAVE
\emph{(\ref{eq:1})} has a unique solution for any $b\in
\mathbb{R}^{n}$ if and only if matrix $(A+B)^{-1}(A-B)$ is a
$P$-matrix.
\end{theorem}
\textbf{Proof.} Since matrix $A+B$ is nonsingular, then the GAVE
(\ref{eq:1}) can be expressed as the following equivalent form
\begin{equation}\label{eq:31}
x^{+}=(A+B)^{-1}(A-B)x^{-}+2(A+B)^{-1}b,
\end{equation}
where  $x^{+}=|x|+x$ and $x^{-}=|x|-x$, which is a  linear
complementarity problem in \cite{Murty}. Thus, according to Lemma
2.1, the linear complementarity problem (\ref{eq:31}) has a unique
solution for any $b\in \mathbb{R}^{n}$. Further, the GAVE
(\ref{eq:1}) for any $b\in \mathbb{R}^{n}$  has a unique solution as
well. $\hfill{} \Box$


Combining Theorem 3.1 with Lemma 2.2, we can get Theorem 3.2.

\begin{theorem}
The GAVE $(\ref{eq:1})$ has a unique solution for any $b\in
\mathbb{R}^{n}$ if and only if matrix $A+B\bar{D}$ is nonsingular
for any diagonal matrix  $\bar{D} = \mbox{diag}(\bar{d}_{i})$ with
$\bar{d}_{i}\in [-1, 1]$.
\end{theorem}
\textbf{Proof.} Since we can express matrix  $\bar{D}$ as
\[
\bar{D}=I-2D,
\]
where  $D=\mbox{diag}(d_{i})$ with $0 \leq d_{i}\leq1$, then matrix
$A+B\bar{D}$ is nonsingular for any diagonal matrix $\bar{D} =
\mbox{diag}(\bar{d}_{i})$ with $\bar{d}_{i}\in [-1, 1]$ (it implies
that $A+B$ is nonsingular), if and only if the  matrix
$(A+B)^{-1}(A+B-2BD)$ with $A,B\in \mathbb{R}^{n\times n}$ is
nonsingular for any diagonal matrix $D=\mbox{diag}(d_{i})$ with $0
\leq d_{i}\leq1$.

By the simply computation, we obtain
\begin{align*}
(A+B)^{-1}(A+B-2BD)&=(A+B)^{-1}(AD-BD+A+B-AD-BD)\\
&=(A+B)^{-1}[(A-B)D+A+B-AD-BD]\\
&=(A+B)^{-1}[(A-B)D+(A+B)(I-D)]\\
&=(A+B)^{-1}(A-B)D+I-D.
\end{align*}
Based on Lemma 2.2, it is easy to know that $(A+B)^{-1}(A-B)$ is a
$P$-matrix. Further, based on Theorem 3.1, the GAVE (\ref{eq:1}) has
a unique solution for any $b\in \mathbb{R}^{n}$. $\hfill{} \Box$

As it is known that all the matrices $A$ and $B$ of the GAVE
$(\ref{eq:1})$ in general are two arbitrary $n \times n$ real
matrices. In this case, if we take $B=0$ in Theorem 3.2, Theorem 3.2
reduces to the fundamental theorem of the linear system $Ax=b$ for
$A \in \mathbb{R}^{n\times n}$: \emph{the linear system $Ax=b$ has a
unique solution for any $b\in \mathbb{R}^{n}$ if and only if matrix
$A \in \mathbb{R}^{n\times n}$ is nonsingular}. In a way, Theorem
3.2 generalizes the necessary and sufficient condition for the
unique solution of the linear system $Ax=b$ with  $A \in
\mathbb{R}^{n\times n}$. Of course, we also know that the linear
system $Ax=b$ has a unique solution for any $b\in \mathbb{C}^{n}$ if
and only if matrix $A \in \mathbb{C}^{n\times n}$ is nonsingular. In
addition, it is noted that Theorem 3.2 implies that matrix $A\in
\mathbb{R}^{n\times n}$ in (\ref{eq:1}) should be nonsingular,
whereas,  matrix $B \in \mathbb{R}^{n\times n}$ in (\ref{eq:1}) is
free. That is to say, matrix  $B \in \mathbb{R}^{n\times n}$ in
(\ref{eq:1}) may be nonsingular or singular.

When $B=I$ in Theorem 3.2, we immediately obtain the necessary and
sufficient condition for the unique solution of the AVE
(\ref{eq:2}), see Corollary 3.1.

\begin{cor}
The AVE $(\ref{eq:2})$ has a unique solution for any $b\in
\mathbb{R}^{n}$ if and only if $A+\bar{D}$ is nonsingular for any
diagonal matrix $\bar{D} = \mbox{diag}(\bar{d}_{i})$ with
$\bar{d}_{i}\in [-1, 1]$.
\end{cor}

By observing the proof of Theorem 3.2, we can obtain the following
necessary and sufficient condition for the unique solution of the
GAVE (\ref{eq:1}) as well.

\begin{cor}
The GAVE $(\ref{eq:1})$ has a unique solution for any $b\in
\mathbb{R}^{n}$ if and only if $A+B-2BD$ is nonsingular for any
diagonal matrix $D=\mbox{diag}(d_{i})$ with $0 \leq d_{i}\leq1$.
\end{cor}

Of course, since matrix $D=\mbox{diag}(d_{i})$ with $0 \leq
d_{i}\leq1$ in Corollary 3.2, then we can use matrix $I-D$ instead
of matrix $D$ in matrix $A+B-2BD$. In this case, we have Corollary
3.3.

\begin{cor}
The GAVE $(\ref{eq:1})$ has a unique solution for any $b\in
\mathbb{R}^{n}$ if and only if $A-B+2BD$ is nonsingular for any
diagonal matrix $D=\mbox{diag}(d_{i})$ with $0 \leq d_{i}\leq1$.
\end{cor}

In  Corollary 3.2 and Corollary 3.3, if we take $B=I$,  then
Corollary 3.2 and Corollary 3.3, respectively, reduces to Theorem
3.3 and Theorem 3.2 in \cite{Wu2}, which are main results in
\cite{Wu2}.

Based on Theorem 3.2,  a series of new sufficient conditions for the
GAVE $(\ref{eq:1})$ and the AVE $(\ref{eq:2})$ can be obtained. Some
of new sufficient conditions for the unique solution of the GAVE
(\ref{eq:1}) or  the AVE $(\ref{eq:2})$ for any $b\in
\mathbb{R}^{n}$  are weaker than the previously published works. For
example, if we express matrix $A+B\bar{D}$ in Theorem 3.2 as
\begin{equation}\label{eq:32}
A+B\bar{D}=A(I+A^{-1}B\bar{D}),
\end{equation}
then Theorem 3.3 can be obtained, which is stated below and its
proof is omitted.

\begin{theorem}
If matrix $A$ in $(\ref{eq:1})$ is nonsingular and satisfies
\begin{equation}\label{eq:33}
\rho(A^{-1}B\bar{D})<1
\end{equation}
for any diagonal matrix  $\bar{D}= \mbox{diag}(\bar{d}_{i})$ with
$\bar{d}_{i}\in [-1, 1]$, then the GAVE $(\ref{eq:1})$ for any $b\in
\mathbb{R}^{n}$ has a unique solution. In addition, if matrix $A$
in $(\ref{eq:2})$ is nonsingular and satisfies
\begin{equation}\label{eq:34}
\rho(A^{-1}\bar{D})<1
\end{equation}
for any diagonal matrix $\bar{D} = \mbox{diag}(\bar{d}_{i})$ with
$\bar{d}_{i}\in [-1, 1]$, then the AVE $(\ref{eq:2})$ for any $b\in
\mathbb{R}^{n}$ has a unique solution.
\end{theorem}

It is easy to know that the condition (\ref{eq:33}) in Theorem 3.3
is slightly weaker than Theorem 2 in \cite{Rohn4}. That is to say,
\[
\rho(A^{-1}B\bar{D})\leq\rho(|A^{-1}B|).
\]
By calculate, we have
\[
A^{-1}B\bar{D}\leq
|A^{-1}B\bar{D}|\leq|A^{-1}B||\bar{D}|\leq|A^{-1}B|.
\]
Further, the condition (\ref{eq:34}) of Theorem 3.3 is is slightly
weaker than $\rho(|A^{-1}|)<1$ in \cite{Rohn4} and
$\|A^{-1}\|_{2}<1$ in \cite{Mangasarian2}.

Using the 2-norm for (\ref{eq:32}), together with the 2-norm of the
matrix equal to the  maximal singular value of the matrix,  the
somewhat stronger sufficient condition is obtained, see Corollary
3.4.

\begin{cor}
If matrix $A$ in  $(\ref{eq:1})$ is nonsingular and satisfies
\begin{equation}\label{eq:35}
\sigma_{1}(A^{-1}B)<1,
\end{equation}
then the GAVE $(\ref{eq:1})$ for any $b\in \mathbb{R}^{n}$  has a
unique solution.
\end{cor}

\textbf{Remark 3.1} When $B=I$,   Corollary 3.4  is  the Proposition
4 in \cite{Mangasarian2} and Corollary 3.1 in \cite{Wu2}, but our
proof is different from the proof in \cite{Mangasarian2} and
\cite{Wu2}.

\textbf{Remark 3.2} Noting that $\sigma_{1}(X)\sigma_{n}(X^{-1})=1$
for non-singular matrix $X$. When $B$ in $(\ref{eq:1})$ is
nonsingular, we can use
\begin{equation}\label{eq:36}
\sigma_{n}(B^{-1}A)>1,
\end{equation}
instead of the condition (\ref{eq:35}).  In this case, the GAVE
$(\ref{eq:1})$ for any $b\in \mathbb{R}^{n}$ has a unique solution
as well.  It is not difficult to find that the condition
(\ref{eq:36}) is slighter weaker than the  condition
\begin{equation}\label{eq:37}
\sigma_{1}(B) < \sigma_{n}(A),
\end{equation}
which was provided in \cite{Wu5}. By the simple computations, we
have
\[
\sigma_{n}(B^{-1}A)\geq\sigma_{n}(B^{-1})\sigma_{n}(A)=\frac{\sigma_{n}(A)}{\sigma_{1}(B)}.
\]
Under the condition (\ref{eq:37}), we get the condition
(\ref{eq:36}). Since the condition (\ref{eq:37}) is slighter weaker
than the following condition
\begin{equation}\label{eq:38}
\sigma_{1}(|B|) < \sigma_{n}(A),
\end{equation}
which was provided in \cite{Rohn5}, it follows that the condition
(\ref{eq:36}) is the weakest condition, compared with the conditions
(\ref{eq:37}) and (\ref{eq:38}). In addition, if $B=I$ in
(\ref{eq:36}), then the condition (\ref{eq:36}) reduces to  the
Proposition 3 (i) in \cite{Mangasarian2} and Theorem 3.6 in
\cite{Wu2}.

It is known that when the smallest singular value $\sigma_{n}$ of
matrix $A$ is greater than $0$, matrix $A$ is nonsingular. Based on
this fact, we have Corollary 3.5.

\begin{cor}
If matrix $A$ in  $(\ref{eq:2})$ satisfies $\sigma_{n}(A+I)>2$, then
the AVE $(\ref{eq:2})$ has a unique solution for any $b\in
\mathbb{R}^{n}$.
\end{cor}
\textbf{Proof.} Based on Corollary 3.2, when the matrix $A+I-2D$ is
nonsingular  for any diagonal matrix $D=\mbox{diag}(d_{i})$ with $0
\leq d_{i}\leq1$, the AVE $(\ref{eq:2})$ has a unique solution for
any $b\in \mathbb{R}^{n}$. Let $\sigma_{n}(A+I-2D)$ indicate the
minimal singular value of the matrix $A+I-2D$. Based on Lemma 2.3,
we have
\[
\sigma_{n}(A+I-2D)\geq\sigma_{n}(A+I)-2\sigma_{1}(D).
\]
Since $\sigma_{1}(D)\leq1$, clearly, when $\sigma_{n}(A+I)>2$, we
have $\sigma_{n}(A+I-2D)>0$. This implies that matrix $A+I-2D$ is
nonsingular. $\hfill{} \Box$

Of course, combining  matrix $A+B-2BD$ ($A-B+2BD$) of Corollary 3.2
(Corollary 3.3) with the approach of the proof of Corollary 3.5,
other sufficient conditions can be obtained as well. Here is
omitted.

\section{Conclusions}

In this paper, we have presented some necessary and sufficient
conditions for the unique solution of the  generalized absolute
value equation (GAVE) $Ax-B|x|=b$ with $A, B\in \mathbb{R}^{n\times
n}$. These results not only address the question of the necessary
and sufficient condition for the unique solution of the GAVE, not
only contain the fundamental theorem for the unique solution of the
linear system $Ax=b$ with $A\in \mathbb{R}^{n\times n}$. Moreover,
some presented new sufficient conditions for the unique solution of
the GAVE are weaker than the previous published works.

\end{document}